\documentclass[a4paper,11pt]{amsart}
\usepackage{graphicx}
\usepackage{mathptmx}
\usepackage{amsmath}
\usepackage{amssymb}
\usepackage{enumitem}
\usepackage{xcolor}

\newmuskip\pFqmuskip

\newcommand*\pFq[6][8]{%
  \begingroup 
  \pFqmuskip=#1mu\relax
  \mathcode`=\string"8000
  \begingroup\lccode`\~=`\,
  \lowercase{\endgroup\let~}\pFqcomma
  F^{#2}_{#3}{\left(\genfrac..{0pt}{}{#4}{#5}\bigg|#6\right)}%
  \endgroup
}
\newcommand{\pFqcomma}{\mskip\pFqmuskip}

\newtheorem{theorem}{Theorem}

\newtheorem{corollary}[theorem]{Corollary}
\newtheorem{proposition}[theorem]{Proposition}

\begin{document}

\title[Complete and incomplete Bell polynomials associated with  Lah-Bell numbers and  polynomials]{Complete and incomplete Bell polynomials associated with  Lah-Bell numbers and  polynomials}

\author{Taekyun  Kim $^{1}$}
\address{$^{1}$ Department of Mathematics, Kwangwoon University, Seoul 139-701, Republic of Korea}
\email{tkkim@kw.ac.kr$^{\dagger}$, luciasconstant@kw.ac.kr}

\author{DAE SAN KIM $^{2}$}
\address{$^{2}$ Department of Mathematics, Sogang University, Seoul 121-742, Republic of Korea}
\email{dskim@sogang.ac.kr}

\author{Lee-Chae  Jang$^{3}$}
\address{$^{3}$ Graduate School of Education, Konkuk University, Seoul 05029, Republic of Korea}
\email{Lcjang@konkuk.ac.kr}

\author{Hyunseok Lee $^{4}$}
\address{$^{4}$ Department of Mathematics, Kwangwoon University, Seoul 139-701, Republic of Korea}
\email{luciasconstant@kw.ac.kr}

\author{ Han-Young Kim $^{5}$}
\address{$^{5}$ Department of Mathematics, Kwangwoon University, Seoul 139-701, Republic of Korea}
\email{gksdud213@kw.ac.kr}


\maketitle

\begin{abstract}
The $n$th $r$-extended Lah-Bell number is defined as the number of ways a set with $n+r$ elements can be partitioned into ordered blocks such that $r$ distinguished elements have to be in distinct ordered blocks. The aim of this paper is to introduce incomplete $r$-extended Lah-Bell polynomials and complete $r$-extended Lah-Bell polynomials respectively as multivariate versions of $r$-Lah numbers and the $r$-extended Lah-Bell numbers and to investigate some properties and identities for these polynomials. From these investigations, we obtain some expressions for the $r$-Lah numbers and the $r$-extended Lah-Bell numbers as finite sums.
\end{abstract}

\section{Introduction}
It is well known that the unsigned Lah number $L(n,k),\ (n \ge k \ge 0)$, counts the number of ways a set with $n$ elements can be partitioned into $k$ non-empty linearly ordered subsets (see [4,7,8]). The $n$-th Lah-Bell number $B_{n}^{L},\ (n\ge 0)$, is the number of ways a set with $n$ elements can be partitioned into non-empty linearly ordered subsets. Thus, we note that
\begin{equation}
B_{n}^{L}=\sum_{k=0}^{n}L(n,k),\quad (n\ge 0),\quad(\mathrm{see}\ [7,8]).\label{1}
\end{equation}
From \eqref{1}, it can be shown that the generating function of Lah-Bell numbers is given by
\begin{equation}
e^{\frac{t}{1-t}}=\sum_{n=0}^{\infty}B_{n}^{L}\frac{t^{n}}{n!},\quad(\mathrm{see}\ [7,8]),\label{2}
\end{equation}
where
\begin{equation}
\frac{1}{k!}\bigg(\frac{t}{1-t}\bigg)^{k}=\sum_{n=k}^{\infty}L(n,k)\frac{t^{n}}{n!},\ (k\ge 0),\quad (\mathrm{see}\ [8,13,15,17]).
\end{equation}  \label{2-1}
\noindent Explicitly, we see from \eqref{2-1} that the Lah numbers are given by
\begin{displaymath}
L(n,k)=\frac{n!}{k!}\binom{n-1}{k-1},\quad(n \ge k \ge 0),\quad(\mathrm{see}\ [7,8,9,10,17,19]).
\end{displaymath} \par
 Let $n,k,r$ be non-negative integers with $n\ge k$. Then, the $r$-Lah number $L_{r}(n,k)$ counts the number of partitions of a set with $n+r$ elements into $k+r$ ordered blocks such that $r$ distinguished elements have to be in distinct ordered blocks (see [17]). The $r$-extended Lah-Bell number $B_{n,r}^{L}$ is defined as the number of ways a set with $n+r$ elements can be partitioned into ordered blocks such that $r$ distinguished elements have to be in distinct ordered blocks (see [7]). By the definitions of $r$-Lah numbers and $r$-extended Lah-Bell numbers, we have
\begin{equation}
B_{n,r}^{L}=\sum_{k=0}^{n}L_{r}(n,k),\quad(n\ge 0),\quad(\mathrm{see}\ [7]).\label{3}
\end{equation}
From \eqref{3}, we see that the generating function of $r$-extended Lah-Bell numbers is given by
\begin{equation}
e^{\frac{t}{1-t}}\bigg(\frac{1}{1-t}\bigg)^{2r}=\sum_{n=0}^{\infty}B_{n,r}^{L}\frac{t^{n}}{n!},\quad(\mathrm{see}\ [7,15]),\label{4}	
\end{equation}
where
\begin{equation}
\frac{1}{k!}\bigg(\frac{t}{1-t}\bigg)^{k}\bigg(\frac{1}{1-t}\bigg)^{2r}=\sum_{n=k}^{\infty}L_{r}(n,k)\frac{t^{n}}{n!},\quad(\mathrm{see}\ [7,17]),\label{5}
\end{equation}
where $k$ is a non-negative integer. \par
\noindent Explicitly, the $r$-Lah numbers are given by
\begin{displaymath}
L_{r}(n,k)=\frac{n!}{k!}\binom{n+2r-1}{k+2r-1},\quad(n \ge k \ge 0),\quad(\mathrm{see}\ [7,8,9,10,17,19]).
\end{displaymath} \par
\indent  In [7], the $r$-extended Lah-Bell polynomials are defined by
\begin{equation}
e^{x(\frac{t}{1-t})}\bigg(\frac{1}{1-t}\bigg)^{2r}=\sum_{n=0}^{\infty}B_{n,r}^{L}(x)\frac{t^{n}}{n!}.\label{6}	
\end{equation}
It is well known that the complete Bell polynomials are defined by
\begin{equation}
\exp\bigg(\sum_{j=1}^{\infty}x_{j}\frac{t^{j}}{j!}\bigg)=\sum_{n=0}^{\infty}B_{n}(x_{1},x_{2},\dots,x_{n})\frac{t^{n}}{n!},\quad(\mathrm{see}\ [2,3,4,5,11,14]).\label{7}	
\end{equation}
Then it can be shown that the complete Bell polynomials are given by
\begin{equation}
B_{n}(x_{1},x_{2},\dots,x_{n})=\sum_{j_{1}+2j_{2}+\cdots+nj_{n}=n}\frac{n!}{j_{1}!j_{2}!\cdots j_{n}!}\bigg(\frac{x_{1}}{1}\bigg)^{j_{1}}\bigg(\frac{x_{2}}{2!}\bigg)^{j_{2}}\cdots\bigg(\frac{x_{n}}{n!}\bigg)^{j_{n}},\label{8}	
\end{equation}
where the sum runs over all nonnegative integers $j_1,j_2,\dots,j_n$ satisfying $j_1+2j_2+\cdots+nj_n=n$.
The incomplete Bell polynomials are given by
\begin{equation}
\frac{1}{k!}\bigg(\sum_{m=1}^{\infty}x_m\frac{t^{m}}{m!}\bigg)^{k}	=\sum_{n=k}^{\infty}B_{n,k}(x_{1},x_{2},\dots,x_{n-k+1})\frac{t^{n}}{n!},\quad(n\ge 0),\quad(\mathrm{see}\ [5,11,14]).\label{9}
\end{equation}
Thus, we note that
\begin{equation}
B_{n,k}(x_{1},x_{2},\dots,x_{n-k+1})=\sum_{\pi(n,k)}\frac{n!}{j_{1}!j_{2}!\cdots j_{n-k+1}!}\bigg(\frac{x_{1}}{1!}\bigg)^{j_{1}}\bigg(\frac{x_{2}}{2!}\bigg)^{j_{2}}\cdots \bigg(\frac{x_{n-k+1}}{(n-k+1)!}\bigg)^{j_{n-k+1}},\label{10}
\end{equation}
where the sum runs over the set $\pi(n,k)$ of all nonnegative integers $(j_{i})_{i\ge 1}$ satisfying $j_{1}+j_{2}+\cdots+j_{n-k+1}=k$, and $1 j_{1}+2 j_{2}+\cdots+(n-k+1)j_{n-k+1}=n$.\par
\noindent Then the complete and incomplete Bell polynomials are related by
\begin{displaymath}
B_{n}(x_{1},x_{2},\dots,x_{n})=\sum_{k=1}^{n}B_{n,k}(x_{1},x_{2},\dots,x_{n-k+1}),\quad(n \ge 1).
\end{displaymath}

Let $f$ be a $C^{\infty}$-function. Then we have
\begin{align}
e^{f(x+t)}\ &=\ \exp\bigg(\sum_{j=0}^{\infty}f^{(j)}(x)\frac{t^{j}}{j!}\bigg)  \label{11}\\
&=\ \exp\bigg(f(x)+\sum_{j=1}^{\infty}f^{(j)}(x)\frac{t^{j}}{j!}\bigg) \nonumber\\
&\ =\ e^{f(x)}\bigg(1+\sum_{n=1}^{\infty}B_{n}\big(f^{(1)}(x),f^{(2)}(x),\dots,f^{(n)}(x)\big)\frac{t^{n}}{n!}\bigg),\nonumber	
\end{align}
where $f^{(j)}(x)$ is the $j$-th derivative of $f(x)$ and $\exp(t)=e^{t}$. \par
We observe that
\begin{align}
\frac{d^{m}}{dx^{m}}e^{f(x)}\ &=\ \frac{\partial^{m}}{\partial x^{m}}e^{f(x+t)}\bigg|_{t=0}\ =\ \frac{\partial^{m}}{\partial t^{m}}e^{f(x+t)}\bigg|_{t=0}\label{12} \\
&=\ e^{f(x)}B_{m}\big(f^{(1)}(x),f^{(2)}(x),\dots,f^{(m)}(x)\big). \nonumber	
\end{align}
From \eqref{11} and \eqref{12}, we obtain K\"olbig and Coeffey's equation as follows :
\begin{equation}
\frac{d^{m}}{dx^{m}}e^{f(x)}=e^{f(x)}B_{m}\big(f^{(1)}(x),f^{(2)}(x),\dots,f^{(m)}(x)\big),\quad(m\ge 1),\quad(\mathrm{see}\ [6,12]).\label{13}
\end{equation}
The exponential incomplete $r$-Bell polynomials are defined by the generating function
\begin{equation}
\frac{1}{k!}\bigg(\sum_{j=1}^{\infty}a_{j}\frac{t^{j}}{j!}\bigg)^{k}\bigg(\sum_{i=0}^{\infty}b_{i+1}\frac{t^{i}}{i!}\bigg)^{r}=\sum_{n=k}^{\infty}B_{n+r,k+r}^{(r)}(a_{1},a_{2},\dots :b_{1},b_{2},\dots)\frac{t^{n}}{n!}. \label{14}	
\end{equation}
From \eqref{14}, we note that
\begin{align}
&\ B_{n+r,k+r}^{(r)}(a_{1},a_{2},\dots:b_{1},b_{2},\dots)\label{15} \\
&\quad=\ \sum_{\Lambda(n,k,r)}\bigg[\frac{n!}{k_{1}!k_{2}!k_{3}\cdots}\bigg(\frac{a_{1}}{1!}\bigg)^{k_{1}} \bigg(\frac{a_{2}}{2!}\bigg)^{k_{2}} \bigg(\frac{a_{3}}{3!}\bigg)^{k_{3}}\cdots\bigg]\nonumber \\
&\qquad\quad\times\bigg[\frac{r!}{r_{0}!r_{1}!r_{2}\cdots}\bigg(\frac{b_{1}}{0!}\bigg)^{r_{0}} \bigg(\frac{b_{2}}{1!}\bigg)^{r_{1}} \bigg(\frac{b_{3}}{2!}\bigg)^{r_{2}}\cdots\bigg],\nonumber	
\end{align}
where $\Lambda(n,k,r)$ denotes the set of all non-negative integers $(k_{i})_{i\ge 1}$ and $(r_{i})_{i\ge 0}$ such that
\begin{displaymath}
	\sum_{i\ge 1}k_{i}=k,\quad\sum_{i\ge 0}r_{i}=r\quad\mathrm{and}\quad \sum_{i\ge 1}i(k_{i}+r_{i})=n,\quad(\mathrm{see}\ [4,5,14]).
\end{displaymath}
Assume that $(a_{i})_{i\ge 1}$ and $(b_{i})_{i\ge 1}$ are sequences of positive integers.\par Then the number $B_{n+r,k+r}^{(r)}(a_{1},a_{2},\dots;b_{1},b_{2},\dots)$ counts the number of partitions of an $(n+r)$-set into $(k+r)$ blocks satisfying:
\begin{itemize}
	\item The first $r$ elements belong to different blocks,
	\item Any block of size $i$ containing no elements from the first $r$ elements can be colored with $a_{i}$ colors,
	\item Any block of size $i$ containing one element from the first $r$ elements can be colored with $b_{i}$ colors.
\end{itemize}
The complete $r$-Bell polynomials are given by
\begin{equation}
\exp\bigg(\sum_{i=1}^{\infty}a_{i}\frac{t^{i}}{i!}\bigg)\bigg(\sum_{j=0}^{\infty}b_{j+1}\frac{t^{j}}{j!}\bigg)^{r}=\sum_{n=0}^{\infty}B_{n}^{(r)}(a_{1},a_{2},\dots;b_{1},b_{2},\dots)\frac{t^{n}}{n!},\label{16}
\end{equation}
(see [4,5,11,14]). \par
By \eqref{15} and \eqref{16}, we get
\begin{equation}
B_{n}^{(r)}(a_{1},a_{2},\dots:b_{1},b_{2},\dots)=\sum_{k=0}^{n}B_{n+r,k+r}^{(r)}(a_{1},a_{2},\dots:b_{1},b_{2},\dots),\quad(\mathrm{see}\ [5]).\label{17}
\end{equation}

The incomplete and complete Bell polynomials have applications to such diverse areas as combinatorics, probability, algebra and analysis. The number of monomials appearing in the incomplete Bell polynomial $B_{n,k}(x_1,x_2,\cdots,x_{n-k+1})$ is the number of partitioning $n$ into $k$ parts and the coefficient of each monomial is the number of partitioning $n$ as the corresponding $k$ parts. Also, the incomplete Bell polynomials $B_{n,k}(x_1,x_2,\cdots,x_{n-k+1})$ appear in the  Fa\`a di Bruno formula concerning higher-order derivatives of composite functions (see [5]). In addition, the incomplete Bell polynomials can be used in constructing sequences of binomial type (see [16]) and there are certain connections between incomplete Bell polynomials and combinatorial Hopf algebras such as the Hopf algebra of word symmetric functions, the Hopf algebra of symmetric functions and the Faà di Bruno algebra, etc (see [1]).
The complete Bell polynomials $B_n(x_1,x_2,\cdots,x_n)$  have applications to probability theory (see [5,12,18]).
Indeed, the $n$th moment $\mu_n=E[X^n]$ of the random variable $X$ is the $n$th complete Bell polynomial in the first $n$ cumulants $\mu_n=B_n(\kappa_1,\kappa_2,\cdots,\kappa_{n})$. The reader can refer to the Ph. D. thesis of Port [18] for many applications to probability theory and combinatorics.
Many special numbers, like Stirling numbers of both kinds, Lah numbers and idempotent numbers, appear in many combinatorial and number theoretic identities involving complete and incomplete Bell polynomials.
We let the reader refer to the Introduction in [11] for further details on these. \par
The incomplete Lah-Bell polynomials (see \eqref{21}) and the complete Lah-Bell polynomials (see \eqref{24}) are respectively multivariate versions of the unsigned Lah numbers and the Lah-Bell numbers.
We note here that the incomplete Bell polynomials (see \eqref{9}) and the incomplete Lah-Bell polynomials are related as given in \eqref{22}, while the complete Bell polynomials (see \eqref{7}) and the complete Lah-Bell polynomials are related as given in \eqref{25}.
The incomplete $r$-extended Lah-Bell polynomials (see \eqref{29}) and the complete $r$-extended Lah-Bell polynomials (see \eqref{31}) are respectively extended versions of the incomplete Lah-Bell polynomials and the complete Lah-Bell polynomials. Further, they are respectively multivariate versions of the $r$-Lah numbers and the $r$-extended Lah-Bell numbers. \par
The aim of this paper is to introduce the incomplete $r$-extended Lah-Bell polynomials and the complete $r$-extended Lah-Bell polynomials and to investigate some properties and identities for these polynomials. From these investigations, we obtain some expressions for the $r$-Lah numbers and the $r$-extended Lah-Bell numbers as finite sums.

\section{Complete and incomplete $r$-extended Lah-Bell polynomials}
Let $f(t)=\frac{t}{1-t}$. Then we have
\begin{equation}
f^{(n)}(t)=\frac{d^{n}}{dt^{n}}f(t)=\frac{n!}{(1-t)^{n+1}},\quad(n\ge 1).\label{18}	
\end{equation}
By \eqref{13}, we get
\begin{equation}
\frac{d^{n}}{dt^{n}}e^{\frac{t}{1-t}}\bigg|_{t=0}=B_{n}\big(1!,2!,\dots,n!\big).	\label{19}
\end{equation}
From \eqref{2}, we note that
\begin{equation}
\frac{d^{n}}{dt^{n}}e^{\frac{t}{1-t}}\bigg|_{t=0}=\frac{d^{n}}{dt^{n}}\sum_{k=0}^{\infty}B_{k}^{L}\frac{t^{k}}{k!}\bigg|_{t=0}=B_{n}^{L}.\label{20}
\end{equation}
Therefore, by \eqref{19} and \eqref{20}, we obtain the following theorem.
\begin{theorem}
For $n\ge 1$, we have
\begin{displaymath}
B_{n}^{L}=B_{n}(1!,2!,\dots,n!)=\sum_{k_{1}+2k_{2}+\cdots+nk_{n}=n}\frac{n!}{k_{1}!k_{2}!\cdots k_{n}!}.
\end{displaymath}
\end{theorem}
Let us consider the incomplete Lah-Bell polynomials given by
\begin{equation}
\frac{1}{k!}\bigg(\sum_{m=1}^{\infty}x_{m}t^{m}\bigg)^{k}=\sum_{n=k}^{\infty}B_{n,k}^{L}\big(x_{1},x_{2},\dots,x_{n-k+1}\big)\frac{t^{n}}{n!}, \label{21}
\end{equation}
where $n,k\ge 0$ with $n\ge k$. \par
Note hat $B_{n,k}^{L}(1,1,\dots,1)=L(n,k)$, $(n \ge k \ge 0)$.\par
Indeed, by \eqref{9} and \eqref{21} , we get
\begin{equation}
B_{n,k}^{L}\big(x_{1},x_{2},\dots,x_{n-k+1}\big)=B_{n,k}\big(1!x_{1},2!x_{2},\dots,(n-k+1)!x_{n-k+1}\big).\label{22}	
\end{equation}
From \eqref{22}, we note that
\begin{align*}
	L(n,k)\ &=\ B_{n,k}^{L}(1,1,\dots,1)\ =\ B_{n,k}\big(1!,2!,\dots,(n-k+1)!\big) \\
	&=\ \sum_{\substack{j_{1}+j_{2}+\cdots+j_{n-k+1}=k\\ j_{1}+2j_{2}+\cdots+(n-k+1)j_{n-k+1}=n}}\frac{n!}{j_{1}!j_{2}!\cdots j_{n-k+1}!}
\end{align*}
Therefore, by \eqref{22}, we obtain the following proposition.
\begin{proposition}
	For $n,k\ge 0$ with $n\ge k$, we have
	\begin{displaymath}
		B_{n,k}^{L}\big(x_{1},x_{2},\dots,x_{n-k+1}\big)=B_{n,k}\big(1!x_{1},2!x_{2},\dots,(n-k+1)x_{n-k+1}\big).
	\end{displaymath}
In addition,
\begin{displaymath}
	L(n,k)= \sum_{\substack{j_{1}+j_{2}+\cdots+j_{n-k+1}=k\\ j_{1}+2j_{2}+\cdots+(n-k+1)j_{n-k+1}=n}}\frac{n!}{j_{1}!j_{2}!\cdots j_{n-k+1}!}.
\end{displaymath}
\end{proposition}
From \eqref{22}, we note that
\begin{align}
B_{n,k}^{L}\big(\alpha x_{1},\alpha x_{2},\dots,\alpha x_{n-k+1}\big)\ &=\ B_{n,k}\big(\alpha 1!x_{1}, \alpha 2!x_{2},\dots,\alpha(n-k+1)!x_{n-k+1}\big)\label{23}\\
&=\ 	\alpha^{k}B_{n,k}\big(1!x_{1}, 2!x_{2},\dots,(n-k+1)!x_{n-k+1}\big)\nonumber \\
&=\ \alpha^{k}B_{n,k}^{L}(x_{1},x_{2},\dots,x_{n-k+1}).\nonumber
\end{align}
We now consider the complete Lah-Bell polynomials given by
\begin{equation}
\exp\bigg(\sum_{i=1}^{\infty}x_{i}t^{i}\bigg)=\sum_{n=0}^{\infty}B_{n}^{L}(x_{1},x_{2},\dots,x_{n})\frac{t^{n}}{n!}.\label{24}	
\end{equation}
By \eqref{24}, we get
\begin{align}
B_{n}^{L}\big(x_{1},x_{2},\dots,x_{n}\big)\ &=\ B_{n}\big(1!x_{1},2!x_{2},\dots,n!x_{n}\big) \label{25}\\
&=\ \sum_{l_{1}+2l_{2}+\cdots+nl_{n-1}=n}\frac{n!}{l_{1}!l_{2}!\cdots l_{n}!}x_{1}^{l_{1}}x_{2}^{l_{2}}\cdots x_{n}^{l_{n}},\quad (n\ge 0).\nonumber
\end{align}
From \eqref{21} and \eqref{24}, we note that
\begin{align}
1+\sum_{n=1}^{\infty}B_{n}^{L}(x_{1},x_{2},\dots,x_{n})\frac{t^n}{n!}&=\exp\bigg(\sum_{i=1}^{\infty}x_{i}t^{i}\bigg) \label{26} \\
&=\ 1+\sum_{k=1}^{\infty}\frac{1}{k!}\bigg(\sum_{i=1}^{\infty}x_{i}t^{i}\bigg)^{k} \nonumber \\
&=\ 1+\sum_{k=1}^{\infty}\sum_{n=k}^{\infty}B_{n}^{L}(x_{1},x_{2},\dots,x_{n-k+1})\frac{t^{n}}{n!}\nonumber \\
&=\ 1+\sum_{n=1}^{\infty}\bigg(\sum_{k=1}^{n}B_{n,k}^{L}(x_{1},x_{2},\dots,x_{n-k+1}\bigg)\frac{t^{n}}{n!}.\nonumber
\end{align}
Therefore, by \eqref{24} and \eqref{26}, we obtain the following theorem.
\begin{theorem}
For $n\ge 1$, we have
\begin{align*}
B_{n}^{L}(x_{1},x_{2},\dots,x_{n})\ &=\ \sum_{k=1}^{n}B_{n,k}^{L}\big(x_{1},x_{2},\dots,x_{n-k+1}\big)\\
&=\ \sum_{k=1}^{n}B_{n,k}\big(1!x_{1},2!x_{2},\dots,(n-k+1)!x_{n-k+1}\big).
\end{align*}
In addition, for $n\ge 1$, we have
\begin{displaymath}
B_{n}^{L}(x_{1},x_{2},\dots,x_{n})=\sum_{k=1}^{n}\sum_{\pi(n,k)}\frac{n!}{l_{1}!l_{2}!\cdots l_{n-k+1}!}x_{1}^{l_{1}}x_{2}^{l_{2}}\cdots x_{n-k+1}^{l_{n-k+1}},
\end{displaymath}
where $\pi(n,k)$ denotes the set of all non-negative integers $(l_{i})_{i\ge 1}$ such that $l_{1}+l_{2}+\cdots+l_{n-k+1}=k$ and $1\cdot l_{1}+2\cdot l_{2}+\cdots+(n-k+1)l_{n-k+1}=n$.
\end{theorem}
By \eqref{24}, we easily get
\begin{equation}
\sum_{n=0}^{\infty}B_{n}^{L}(1,1,\dots,1)\frac{t^{n}}{n!}=\exp\bigg(\sum_{i=1}^{\infty}t^{i}\bigg)=\exp\bigg(\frac{t}{1-t}\bigg)=\sum_{n=0}^{\infty}B_{n}^{L}\frac{t^{n}}{n!}.\label{27}
\end{equation}
From \eqref{27}, we note that
\begin{displaymath}
B_{n}^{L}(1,1,\dots,1)=B_{n}^{L},\quad (n\ge 0).
\end{displaymath}
By Proposition 2, \eqref{23} and Theorem 3, we get
\begin{align}
B_{n}^{L}(x,x,\dots,x)\ &=\ \sum_{k=0}^{n}B_{n,k}^{L}(x,x,\dots,x)\ =\ \sum_{k=0}^{n}B_{n,k}\big(1!x,2!x,\dots,(n-k+1)!x_{n-k+1}\big) \label{28} \\
&=\ \sum_{k=0}^{n}x^{k}B_{n,k}\big(1!,2!,\dots,(n-k+1)!\big)	\ =\ \sum_{k=0}^{n}x^{k}L(n,k)\ =\ B_{n}^{L}(x).\nonumber
\end{align}
Assume that $\{a_{i}\}_{i\ge 1}$ and $\{b_{i}\}_{i\ge 1}$ are sequences of positive integers. We define the {\it{incomplete $r$-extended Lah-Bell polynomials}} by
\begin{equation}
\frac{1}{k!}\bigg(\sum_{j=1}^{\infty}a_{j}t^{j}\bigg)^{k}\bigg(\sum_{i=0}^{\infty}b_{i+1}t^{i}\bigg)^{2r}=\sum_{n=k}^{\infty}B_{n+2r,k+2r}^{L}(a_{1},a_{2},\dots:b_{1},b_{2},\dots)\frac{t^{n}}{n!},\label{29}
\end{equation}
where $k,r$ are nonnegative integers. \par
From \eqref{29}, we have
\begin{align}
& B_{n+2r,k+2r}^{L}(a_{1},a_{2},\dots:b_{1},b_{2},\dots)=B_{n+2r,k+2r}^{(2r)}(1!a_{1},2!a_{2},\dots:0!b_{1},1!b_{2},\dots) \label{30}\\
&\ =\sum_{\Lambda(n,k,2r)}\bigg[\frac{n!}{k_{1}!k_{2}!k_{3}!\cdots}a_{1}^{k_{1}}a_{2}^{k_{2}}a_{3}^{k_{3}}\cdots\bigg]\bigg[\frac{(2r)!}{r_{0}!r_{1}!r_{2}\cdots}b_{1}^{r_{0}}b_{2}^{r_{1}}b_{3}^{r_{2}}\cdots \bigg],\nonumber
\end{align}
where $\Lambda(n,k,2r)$ denotes the set of all nonnegative integers $\{k_{i}\}_{i\ge 1}$ and $\{r_{i}\}_{i\ge 0}$ such that $\displaystyle\sum_{i\ge 1}k_{i}=k,\ \sum_{i\ge 0}r_{i}=2r\displaystyle$ and $\displaystyle\sum_{i\ge 1}i(k_{i}+r_{i})=n\displaystyle$. \par
We define  the {\it{complete $r$-extended Lah-Bell polynomials}} $B_{n}^{(L,2r)}(x|a_{1},a_{2},\dots:b_{1},b_{2},\dots)$, $(n\ge 0)$, which are given by
\begin{equation}
\exp\bigg(x\sum_{j=1}^{\infty}a_{j}t^{j}\bigg)\bigg(\sum_{i=0}^{\infty}b_{i+1}t^{i}\bigg)^{2r}=\sum_{n=0}^{\infty}B_{n}^{(L,2r)}(x|a_{1},a_{2},\dots:b_{1},b_{2},\dots)\frac{t^{n}}{n!}. \label{31}
\end{equation}
Thus, we note that
\begin{align}
& \exp\bigg(x\sum_{j=1}^{\infty}a_{j}t^{j}\bigg)\bigg(\sum_{i=0}^{\infty}b_{i+1}t^{i}\bigg)^{2r}\ =\ \sum_{k=0}^{\infty}\frac{x^{k}}{k!}\bigg(\sum_{j=1}^{\infty}a_{j}t^{j}\bigg)^{k}\bigg(\sum_{i=0}^{\infty}b_{i+1}t^{i}\bigg)^{2r}\label{32} \\
&\ = \sum_{k=0}^{\infty}x^{k}\sum_{n=k}^{\infty}B_{n+2r,k+2r}^{L}(a_{1},a_{2},\dots:b_{1},b_{2},\dots)\frac{t^{n}}{n!}\ =\ \sum_{n=0}^{\infty}\sum_{k=0}^{n}x^{k}B_{n+2r,k+2r}^{L}(a_{1},a_{2},\dots:b_{1},b_{2},\dots)\frac{t^{n}}{n!}.\nonumber	
\end{align}
 From \eqref{31} and \eqref{32}, we have
\begin{equation}
B_{n}^{(L,2r)}\big(x\ |\ a_{1},a_{2},\dots:b_{1},b_{2},\dots)=\sum_{k=0}^{n}x^{k}B_{n+2r,k+2r}^{L}(a_{1},a_{2},\dots:b_{1},b_{2},\dots),\ (n\ge 0). \label{33}
\end{equation}
By \eqref{17}, \eqref{30}, \eqref{31} and \eqref{33}, we have
\begin{align}
&B_{n}^{(2r)}(1!a_{1},2!a_{2},\dots:0!b_{1},1!b_{2},\dots)
=\sum_{k=0}^{n}B_{n+2r,k+2r}^{(2r)}(1!a_{1},2!a_{2},\dots:0!b_{1},1!b_{2},\dots) \label{34}\\
&= \sum_{k=0}^{n}B_{n+2r,k+2r}^{L}(a_{1},a_{2},\dots:b_{1},b_{2},\dots)
= B_{n}^{(L,2r)}(1|a_{1},a_{2},\dots:b_{1},b_{2},\dots).\nonumber
\end{align}
Therefore, by \eqref{30} and \eqref{33}, we obtain the following theorem.
\begin{theorem}
For $n\ge 0$, we have
\begin{align*}
B_{n}^{(L,2r)}\big(x\ |\ a_{1},a_{2},\dots:b_{1},b_{2},\dots\big)\ &=\ \sum_{k=0}^{n}x^{k}B_{n+2r,k+2r}^{L}(a_{1},a_{2},\dots:b_{1},b_{2},\dots)\\
&=\ \sum_{k=0}^{n}x^{k}B_{n+2r,k+2r}^{(2r)}(1!a_{1},2!a_{2},\dots:0!b_{1},1!b_{2},\dots).
\end{align*}
\end{theorem}
From \eqref{29}, we note that
\begin{align}
&\sum_{n=k}^{\infty}B_{n+2r,k+2r}^{L}\big(1,1,\dots;1,1,\dots\big)\frac{t^{n}}{n!}\ =\ \frac{1}{k!}\bigg(\sum_{j=1}^{\infty}t^{j}\bigg)^{k}\bigg(\sum_{i=0}^{\infty}t^{i}\bigg)^{2r}\label{36}\\
&\ =\ \frac{1}{k!}	\bigg(\frac{t}{1-t}\bigg)^{k}\bigg(\frac{1}{1-t}\bigg)^{2r}\ =\ \sum_{n=k}^{\infty}L_{r}(n,k)\frac{t^{n}}{n!}.\nonumber
\end{align}
By \eqref{31} and \eqref{36}, we get
\begin{align}
&\sum_{n=0}^{\infty}B_{n}^{(L,2r)}(x|1,1,\dots:1,1,\dots)\frac{t^{n}}{n!}\ =\ \exp\bigg(x\sum_{j=1}^{\infty}t^{j}\bigg)\bigg(\sum_{i=0}^{\infty}t^{i}\bigg)^{2r}\label{37} \\
&=\ e^{x(\frac{t}{1-t})}\cdot\bigg(\frac{1}{1-t}\bigg)^{2r}\ =\ \sum_{n=0}^{\infty}\bigg(\sum_{k=0}^{n}x^{k}L_{r}(n,k)\bigg)\frac{t^{n}}{n!}. \nonumber
\end{align}
Thus, by \eqref{36} and \eqref{37}, we have
\begin{align}
B_{n}^{(L,2r)}\big(x\ |\ 1,1,\dots:1,1,\dots)\ &=\ \sum_{k=0}^{n}x^{k}B_{n+2r,k+2r}^{L}(1,1,\dots:1,1,\dots)\label{38}\\
&=\ \sum_{k=0}^{n}x^{k}L_{r}(n,k).\nonumber
\end{align}
Therefore, we obtain the following theorem.
\begin{theorem}
For $n \ge k \ge 0$, we have
\begin{displaymath}
B_{n+2r,k+2r}^{L}(1,1,\dots:1,1,\dots)=L_{r}(n,k)
\end{displaymath}
and
\begin{displaymath}
B_{n}^{(L,2r)}(x\ |\ 1,1,\dots:1,1,\dots)=\sum_{k=0}^{n}B_{n+2r,k+2r}^{L}(1,1,\dots:1,1,\dots)=\sum_{k=0}^{n}x^{k}L_{r}(n,k).
\end{displaymath}
\end{theorem}
From \eqref{36} and \eqref{30}, we note that
\begin{align*}
	L_{r}(n,k)\ &=\ B_{n+2r,k+2r}^{L}(1,1,\dots:1,1,\dots)\\
	&=\ B_{n+2r,k+2r}^{(2r)}(1!,2!,\dots:1!,2!,\dots) \\
	&=\ \sum_{\Lambda(n,k,2r)}\frac{n!}{k_{1}!k_{2}!\cdots}\frac{(2r)!}{r_{0}!r_{1}!\cdots}.
\end{align*}
\begin{corollary}
For $n,k,r\ge 0$ with $n \ge k$, we have
\begin{displaymath}
L_{r}(n,k)= \sum_{\Lambda(n,k,2r)}\frac{n!}{k_{1}!k_{2}!\cdots}\cdot\frac{(2r)!}{r_{0}!r_{1}!\cdots},
\end{displaymath}
where $\Lambda(n,k,2r)$ denotes the set of all nonnegative integers $\{k_{i}\}_{i\ge 1}$ and $\{r_{i}\}_{i\ge }$ such that \par
\noindent $\displaystyle\sum_{i\ge 1}k_{i}=k,\ \sum_{i\ge 0}r_{i}=2r\displaystyle$ and $\displaystyle\sum_{i\ge 1}i(k_{i}+r_{i})=n.\displaystyle$
\end{corollary}

Now, we observe that
\begin{align}
\exp\bigg(\sum_{i=1}^{\infty}x_{i}t^{i}\bigg)\ &=\ 1+\sum_{k=1}^{\infty}\frac{1}{k!}\bigg(\sum_{i=1}^{\infty}x_{i}t^{i}\bigg)^{k} \label{39}	 \\
&=\ 1+\frac{1}{1!}\sum_{i=1}^{\infty}x_{i}t^{i}+\frac{1}{2!}\bigg(\sum_{i=1}^{\infty}x_{i}t^{i}\bigg)^{2}+\frac{1}{3!}\bigg(\sum_{i=1}^{\infty}x_{i}t^{i}\bigg)^{3}+\cdots\nonumber \\
&=\ \sum_{k=0}^{\infty}\sum_{m_{1}+2m_{2}+\cdots+km_{k}=k}\frac{1}{m_{1}!m_{2}!\cdots m_{k}!}x_{1}^{m_{1}}x_{2}^{m_{2}}\cdots x_{k}^{m_{k}}t^{k},\nonumber
\end{align}
and
\begin{equation}	\bigg(\sum_{j=0}^{\infty}y_{j+1}t^{j}\bigg)^{2r}=\sum_{m=0}^{\infty}\sum_{l_{1}+\cdots+l_{2r}=m}
y_{l_{1}+1}y_{l_{2}+1}\cdots y_{l_{2r}+1}t^{m}.\label{40}
\end{equation}
By \eqref{39} and \eqref{40}, we get
\begin{align}
&\exp\bigg(\sum_{i=1}^{\infty}x_{i}t^{i}\bigg)\bigg(\sum_{j=0}^{\infty}y_{j+1}\frac{t^{j}}{j!}\bigg)^{2r}\label{41}\\
&=\ \sum_{n=0}^{\infty}n!\bigg(\sum_{k=0}^{n}\sum_{m_{1}+2m_{2}+\cdots+km_{k}=k}\,\sum_{l_{1}+l_{2}+\cdots+l_{2r}=n-k}\frac{1}{m_{1}!m_{2}!\cdots m_{k}!}x_{1}^{m_{1}}x_{2}^{m_{2}}\cdots x_{k}^{m_{k}}\nonumber \\
&\qquad\quad\times y_{l_{1}+1}y_{l_{2}+1}\cdots y_{l_{2r}+1}\bigg)\frac{t^{n}}{n!}. \nonumber	
\end{align}
Therefore, by \eqref{31} and \eqref{41}, we obtain the following theorem.
\begin{theorem}
For $n,r\ge 0$, we have
\begin{align*}
&B_{n}^{(L,2r)}(1\ |\ x_{1},x_{2},\dots:y_{1},y_{2},\dots)\\
&=\ n!\sum_{k=0}^{n}\sum_{m_{1}+2m_{2}+\cdots+km_{k}=k}\,\sum_{l_{1}+l_{2}+\cdots+l_{2r}=n-k}\frac{1}{m_{1}!m_{2}!\cdots m_{k}!}x_{1}^{m_{1}}\cdots x_{k}^{m_{k}}y_{l_{1}+1}y_{l_{2}+1}\cdots y_{l_{2r}+1}.
\end{align*}
\end{theorem}
\noindent\emph{Remark.} For $n\ge 0$, we have
\begin{equation}
B_{n+2r,k+2r}^{L}(x,x,\dots:1,1,\dots)=x^{k}B_{n+2r,k+2r}^{L}(1,1,\dots:1,1,\dots).\label{42}
\end{equation}
Thus, we note that
\begin{displaymath}
\sum_{k=0}^{n}B_{n+2r,k+2r}^{L}(x,x,\dots:1,1,\dots)=B_{n,r}^{L}(x),\quad(n\ge 0).
\end{displaymath}

\section{Conclusion}
There are various methods of studying special numbers and polynomials, for example,
generating functions, combinatorial methods, umbral calculus, $p$-adic analysis, differential
equations, probability theory, orthogonal polynomials, and special functions. These ways
of investigating special polynomials and numbers can be applied
also to degenerate versions of such polynomials and numbers. Indeed, in recent years,
many mathematicians have drawn their attention to studies of degenerate versions of many
special polynomials and numbers by making use of the aforementioned means
[9,10,14 and references therein]. \par
The incomplete and complete Bell polynomials arise in many different contexts as we stated in the Introduction. For instance, many special numbers, like Stirling numbers of both kinds, Lah numbers and idempotent numbers, appear in many combinatorial and number theoretic identities involving complete and incomplete Bell polynomials. \par
In this paper, we introduced the incomplete $r$-extended Lah-Bell polynomials and the complete $r$-extended Lah-Bell polynomials respectively as multivariate versions of $r$-Lah numbers and the $r$-extended Lah-Bell numbers and investigated some properties and identities for these polynomials. As corollaries to these results, we obtained some expressions for the $r$-Lah numbers and the $r$-extended Lah-Bell numbers as finite sums. \par
It would be very interesting to explore many applications of the incomplete and complete $r$-extended Lah-Bell polynomials just as the incomplete and complete Bell polynomials have diverse applications.

\vspace{0.2in}
{\bf Acknowledgements}

The authors would like to thank the reviewers for their valuable comments and suggestions and Jangjeon Institute for
Mathematical Science for the support of this research.

\vspace{0.2in}
{\bf Funding} 

Not applicable.

\vspace{0.1in}

{\bf Availability of data and materials}

Not applicable.

\vspace{0.1in}

{\bf Competing interests}

The authors declare that they have no conflicts of interest.

\vspace{0.1in}

{\bf Authors’ contributions}

TK and DSK conceived of the framework and structured the whole paper; DSK and TK wrote the paper; LCJ, HL, and HYK
checked the results of the paper; DSK and TK completed the revision of the paper. All authors have read and approved
the final version of the manuscript.

\end{document}